\documentclass[10pt]{article}
\usepackage{amsmath,amssymb,amsfonts,amsthm,mathtools,epsfig}
\usepackage{tikz}
\usepackage{caption,subcaption,verbatim}

\usepackage[bottom=1.5in]{geometry}

\theoremstyle{definition}

\date{ }

\title{\bf {Supercharacter Table of Certain Finite Groups} }

\date{}
\begin{document}

\maketitle

\begin{center} 
\textbf{H. Saydi$^{1,a}$ , M. R. Darafsheh$^{2,b}$ and A. Iranmanesh$^{1,c}$}
\end{center}
\textit{\begin{center}
(1) College of Mathematical Science, Tarbiat Modares University, Tehran, Iran
\\
(2) School of mathematics, statistics and computer science, College of science, University of Tehran, Tehran, Iran
\\
(a) e-mail: h.seydi@modares.ac.ir
\\
(b) e-mail: darafsheh@ut.ac.ir
\\
(c) e-mail: iranmana@modares.ac.ir
\end{center}}

\begin{abstract}
\vskip 3mm
Supercharacter theory is developed by P. Diaconis and I. M. Isaacs as a natural generalization of the classical ordinary character theory. Some classical sums of number theory appear as supercharacters which are obtained by the action of certain subgroups of $GL_d(\mathbb{Z}_n)$ on $\mathbb{Z}_n^d$. In this paper we take $\mathbb{Z}_p^d$, $p$ prime, and by the action of certain subgroups of $GL_d(\mathbb{Z}_p)$ we find supercharacter table of $\mathbb{Z}_p^d$. 

\end{abstract}

\section{Introduction}
 Let $Irr(G)$ denote the set of all the irreducible complex characters of a finite group $G$, and let $Con(G)$ denote the set of all the conjugacy classes of $G$. The identity element of $G$ is denoted by $1$ and the trivial character is denoted by $1_G$. By definition a supercharacter theory for $G$ is a pair $(\mathcal{X},\mathcal{K})$ where $\mathcal{X}$ and $\mathcal{K}$ are partitions of $Irr(G)$ and $G$ respectively, $\vert\mathcal{X}\vert=\vert\mathcal{K}\vert$, $\lbrace1\rbrace\in\mathcal{K}$, and for each $X\in\mathcal{X}$ there is a character $\sigma_X$ such that $\sigma_X(x)=\sigma_X(y)$ for all $x,y\in K$, $K\in\mathcal{K}$. We call $\sigma_X$ supercharacter and each member of $\mathcal{K}$ a superclass. We write $Sup(G)$ for the set of all the supercharacter theories of $G$. 
 
Supercharacter theory of a finite group were defined by Diaconis and Isaacs \cite{DI} as a general case of the ordinary character theory. In fact, in a supercharacter theory characters play the role of irreducible ordinary characters and union of conjugacy classes play the role of conjugacy classes. In \cite{DI} it is shown that $\lbrace1_G\rbrace\in\mathcal{X}$ and if $X\in\mathcal{X}$ then $\sigma_X$ is a constant multiple of $\underset{\chi\in X}{\sum}\chi(1)\chi$, and that we may assume that $\sigma_X=\underset{\chi\in X}{\sum}\chi(1)\chi$. 

For any finite group there are two trivial supercharacter theories as follows. In the first case, $\mathcal{X}=\underset{\chi\in Irr(G)}{\bigcup\lbrace\chi\rbrace}$ and $\mathcal{K}$ is the set of all conjugacy classes of $G$. In the second case, 
\begin{center}
$\mathcal{X}=\lbrace1_G\rbrace\cup\lbrace Irr(G)-\lbrace1_G\rbrace\rbrace$
\end{center}
 and $\mathcal{K}=\lbrace1\rbrace\cup\lbrace G-\lbrace1\rbrace\rbrace$. In the first case, supercharacters are just irreducible characters and superclasses are conjugacy classes. In the second case, the non-trivial supercharacter is $\rho_G-1_G$, where $\rho_G$ denotes the regular character of $G$. These two supercharacter theories of $G$ are denoted by $m(G)$ and $M(G)$ respectively. 

It is mentioned in \cite{HE} that the set of supercharacter theories of a group form a Lattice in the following natural way. $Sup(G)$ can be made to a poset by defining $(\mathcal{X},\mathcal{K})\leq(\mathcal{Y},\mathcal{L})$ if $\mathcal{X}\leq\mathcal{Y}$ in the sense that every part of $\mathcal{X}$ is a subset of some part of $\mathcal{Y}$. In \cite{HE} it is shown that this definition is equivalent to $(\mathcal{X},\mathcal{K})\leq(\mathcal{Y},\mathcal{L})$ if $\mathcal{K}\leq\mathcal{L}$. By this definition $m(G)$ is the least and $M(G)$ is the largest element of $Sup(G)$. 

Among construction of supercharacter theories of a finite group $G$ the following is of great importance which is a lemma by Brauer on character tables of groups. Let $A$ be a subgroup of $\mathbb{A}ut(G)$ and $Irr(G)=\lbrace \chi_1=1_G,\ldots,\chi_h\rbrace$, $Con(G)=\lbrace \mathcal{C}_1=\lbrace1\rbrace,\ldots,\mathcal{C}_h\rbrace$. Suppose for each $\alpha\in A$, ${\mathcal{C}_i}^\alpha=\mathcal{C}_j$, $1\leq i\leq h$, and ${\chi_i}^\alpha(g)=\chi_i(g^\alpha)$ for all $g\in G$, $\alpha\in A$, then the number of conjugacy classes fixed by $ \alpha $ equals the number of irreducible characters fixed by $ \alpha $, and more over the number of orbits of $A$ on $Con(G)$ equals the number of orbits of $A$ on $Irr(G)$, \cite{DO}. It is easy to see that the orbits of $A$ on $Irr(G)$ and $Con(G)$ yield a supercharacter theory for $G$. This supercharacter theory of $G$ is called automorphic. In \cite{HEN} it is shown that all the supercharacter theories of the cyclic group of order $p$, $p$ prime, are automorphic.

 Another aspect of the supercharacter theory of finite groups is to employ the theory to the group $U_n(F)$, the group of $n\times n$ unimodular upper triangular matrices over the Galois field $GF(p^m)$, $p$ prime. Computation of the conjugacy classes and irreducible characters of $U_n(F)$ is still open, but in \cite{AN} the author has developed an applicable supercharacter theory for $U_n(F)$. This result is reviewed in \cite{DI}. 
 
\section{Supercharacter table}
 Let $G$ be a finite group and $(\mathcal{X},\mathcal{K})$ be a supercharacter theory for $G$. Suppose $\mathcal{X}=\lbrace X_1,X_2,\ldots,X_h\rbrace$ be a partition for $Irr(G)$ with the corresponding supercharacter $\sigma_i=\underset{\chi\in X_i}{\sum}\chi(1)\chi$. 

Let $\mathcal{K}=\lbrace K_1,K_2,\ldots,K_h\rbrace$ be the partion of $G$ into superclasses. In fact $X_1=\lbrace1_G\rbrace$ , $K_1=\lbrace1\rbrace$ and $K_i$'s are conjugacy classes of $G$. The supercharacter table of $G$ corresponding to $(\mathcal{X},\mathcal{K})$ is the following $h\times h$ array: 
\begin{center}
\begin{tabular}{|c|c|c|c|c|c|c|}\hline
 & $K_1$&$ K_2$&$ \cdots$&$ K_j$&$ \cdots$&$ K_h$  \\ \hline
 $\sigma_1$ & $\sigma_1(K_1)$&$ \sigma_1(K_2)$&$ \cdots$&$ \sigma_1(K_j)$&$ \cdots$&$ \sigma_1(K_h)$  \\ 
 $\sigma_2$ & $\sigma_2(K_1)$&$ \sigma_2(K_2)$&$ \cdots$&$ \sigma_2(K_j)$&$ \cdots$&$ \sigma_2(K_h)$  \\ 
 $\vdots$ & $\vdots$&$ \vdots$&&&&$ \vdots$  \\ 
 $\sigma_i$ & $\sigma_i(K_1)$&$ \sigma_i(K_2)$&$ \cdots$&$ \sigma_i(K_j)$&$ \cdots$&$ \sigma_i(K_h)$  \\ 
 $\vdots$ & $\vdots$&$ \vdots$&&&&$ \vdots$  \\ 
 $\sigma_h$ & $\sigma_h(K_1)$&$ \sigma_h(K_2)$&$ \cdots$&$ \sigma_h(K_j)$&$ \cdots$&$ \sigma_h(K_h)$  \\  \hline
\end{tabular}
\end{center}

Let us set $S=(\sigma_i(K_j))^h_{i,j=1}$, and call it the supercharacter table of $G$. 

Recall that a class function on $G$ is a function $f:G\longrightarrow\mathbb{C}$ which is constant on conjugacy classes of $G$. The set of all the class functions on $G$, $Cf(G)$ has the structure of a vector space over $\mathbb{C}$ with an orthonormal basis $Irr(G)$ with respect the inner product $\langle f,g\rangle=\dfrac{1}{\vert G\vert}\underset{x\in G}{\sum}f(x)\overline{g(x)}$. Since supercharacters are constant on superclasses, it is natural to call them superclass functions. We have: 
$$ \langle \sigma_i,\sigma_j\rangle=\dfrac{1}{\vert G\vert}\underset{k=1}{\overset{h}{\sum}}\vert K_k\vert\sigma_i(K_k)\overline{\sigma_j(K_k)} $$
But using the orthogonality of $Irr(G)$ we also can write: 
$$ \langle \sigma_i,\sigma_j\rangle=\langle\underset{\chi\in X_i}{\sum}\chi(1)\chi,\underset{\varphi\in X_j}{\sum}\varphi(1)\varphi\rangle=\delta_{ij}\underset{\chi\in X_i}{\sum}\chi(1)^2$$
Therefore 
$$\dfrac{1}{\vert G\vert}\underset{k=1}{\overset{h}{\sum}}\vert K_k\vert\sigma_i(K_k)\overline{\sigma_j(K_k)} =\delta_{ij}\underset{\chi\in X_i}{\sum}\chi(1)^2.$$
If we set the matrix 
$$ U=\dfrac{1}{\sqrt{\vert G\vert}}\begin{bmatrix}\dfrac{\sigma_i(K_j)\sqrt{\vert K_j\vert}}{\sqrt{\underset{\chi\in X_i}{\sum}\chi(1)^2}}\end{bmatrix}^h_{i,j=1} $$
We see that $U$ is a unitary matrix with the following properties, which are proved in \cite{BR}. We have $U=U^t$, $U^2=P$ where $P$ is a permutation matrix and $U^4=I$. 

In the course of studying the supercharacter theory of a group $G$ finding the supercharacter table of $G$ and the matrix $U$ is of great importance. In this paper, we will do this task for certain groups acting on certain sets. 

\section{Automorphic supercharacter table}
 In this section we follow the method used in \cite{BR} considering the group $G=\mathbb{Z}_n^d$ which abelian of order $n^d$. The automophism group of $G$ is $GL_d(\mathbb{Z}_n)$, the group of $d\times d$ invertible matrices with entries in $\mathbb{Z}_n$. We write elements of $G$ as row vectors $y=(y_1,\ldots,y_d)$ and let the action of $GL_d(\mathbb{Z}_n)$ on $G$ be as follows: 
\begin{center}
$y^{A}=yA$
for $A\in GL_d(\mathbb{Z}_n)$.
\end{center}

Irreducible characters of $G$ are of degree $1$ and the number of them is equal to $\vert G\vert$. For $x\in G$, let us define $\psi_x:G\longrightarrow\mathbb{C}^\times$ , by $\psi_x(\zeta)=e(\dfrac{x\cdot \zeta}{n})$, where $e(t)$ stands for $e(t)=e^{2\pi it}$ and $x\cdot \zeta$ is the inner product of two elements $x$ and $\zeta$ of $G$ as row vectors in $G=\mathbb{Z}_n^d$. Therefore $Irr(G)=\lbrace \psi_x\vert x\in G\rbrace$ and the action of $GL_d(\mathbb{Z}_n)$ on $Irr(G)$ is as follows:
 \begin{center}
$\psi_x^A=\psi_{xA^{-t}}$ where $A\in GL_d(\mathbb{Z}_n)$, $x\in G$. 
\end{center}
 
 Now let $\Gamma$ be a symmetric subgroup of $GL_d(\mathbb{Z}_n)$, i. e. $\Gamma^t=\Gamma$. Then $\Gamma$ acts on $G$ and $Irr(G)$ as above Let $\mathcal{X}$ be the set of orbits of $\Gamma$ on $Irr(G)$ and $\mathcal{K}$ be the set of orbits of $GL_d(\mathbb{Z}_n)$ on $Irr(G)$. 
 
It is shown in \cite{BR} that $(\mathcal{X},\mathcal{K})$ is a supercharacter theory of $G$. Following the notations used in \cite{BR} we identity $\psi_x$ with $x$ and $\psi_x^A=\psi_{xA^{-t}}=xA^{-t}$. Therefore $\mathcal{X}$ is identified with the set of orbits of $GL_d(\mathbb{Z}_n)$ on $G$, by $x\longmapsto xA^{-t}$, and $\mathcal{K}$ is identified with the orbits of the action of $GL_d(\mathbb{Z}_n)$ on $G$ by $y\longmapsto yA$. 

In \cite{BR} using different subgroups of $GL_d(\mathbb{Z}_n)$ the authors provide supercharacter tables for $G$. For example the discrete Fourier transform in the case of $ \Gamma=\lbrace1\rbrace $, or $ \Gamma=\lbrace\pm1\rbrace $ a group of order $ 2 $. The Gauss sums is obtained in the case of $ G=\mathbb{Z}_p $, $p$ an odd prime, $\Gamma=\langle g^2\rangle$ where $g$ is a primitive root modulo $p$. Kloosterman sums in the case $ G=\mathbb{Z}_p^2 $, $p$ an odd prime and $\Gamma=\lbrace \begin{bmatrix}
a&0\\0&a^{-1}
\end{bmatrix}\mid 0\neq a\in \mathbb{Z}_p\rbrace$. Heilbronn sums, in the case of $ G=\mathbb{Z}_p^2 $ and $\Gamma=\lbrace x^p\mid 0\neq x\in \mathbb{Z}_p\rbrace$. The Ramanujan sums in the case of $ G=\mathbb{Z}_n $ and $ \Gamma=\mathbb{Z}_n^\times $. It is worth mentioning that all the above sums appear as supercharacters.
 
As a generalization of the group $\Gamma$ in Kloosterman sum we let 
$$ \Gamma=\lbrace \begin{bmatrix}
a&0\\0&b
\end{bmatrix}\mid a,b\in \mathbb{Z}_p^\times\rbrace $$
a group of order $(p-1)^2$. Here $G=\mathbb{Z}_p\times\mathbb{Z}_p$ and orbits of $\Gamma$ on $G$ are: 
\begin{align*}
& Y_1=\lbrace(0,0)\rbrace &  \text{of size}\ & 1\\
& Y_2=(1,0)\Gamma=\lbrace(a,0)\mid a\in \mathbb{Z}_p^\times\rbrace & \text{of size} \ & p-1\\
& Y_3=(0,1)\Gamma=\lbrace(0,b)\mid b\in \mathbb{Z}_p^\times\rbrace & \text{of size} \ & p-1\\
& Y_4=(1,1)\Gamma=\lbrace(a,b)\mid a,b\in \mathbb{Z}_p^\times\rbrace & \text{of size} \ & (p-1)^2
\end{align*}

Orbits of $\Gamma$ on $Irr(G)$ are as follows: 
\begin{align*}
& X_1=\lbrace(0,0)\rbrace &  \text{of size}\ & 1\\
& X_2=(1,0)\Gamma=\lbrace(a,0)\mid a\in \mathbb{Z}_p^\times\rbrace & \text{of size} \ & p-1\\
& X_3=(0,1)\Gamma=\lbrace(0,b)\mid b\in \mathbb{Z}_p^\times\rbrace & \text{of size} \ & p-1\\
& X_4=(1,1)\Gamma=\lbrace(a,b)\mid a,b\in \mathbb{Z}_p^\times\rbrace & \text{of size} \ & (p-1)^2
\end{align*}

Now we form the supercharacter table of $G\cong\mathbb{Z}_p\times\mathbb{Z}_p$. Let $\sigma_i$ be the supercharacter associated with $X_i$, with $\sigma_1=1$. 

We know $\sigma_i=\underset{\psi_{x_i}\in X_i}{\sum}\psi_{x_i}$, and for $y\in Y_j$, $\sigma_i(y)=\underset{\psi_{x_i}\in X_i}{\sum}\psi_{x_i}(y)=\underset{x_i\in X_i}{\sum}e(\dfrac{x_i\cdot y}{p})$, $1\leq i\leq4$. Therefor the following table is calculated: 
\begin{center}
\begin{tabular}{|c|c|c|c|c|}\hline
 $\mathbb{Z}_p^2$& $Y_1$&$ Y_2$&$Y_3$&$Y_4$ \\ 
  superclass size & $1$&$p-1$&$p-1$&$(p-1)^2$ \\ \hline
 $\sigma_1$& $1$&$ 1$&$1$&$1$ \\
 $\sigma_2$& $p-1$&$ -1$&$p-1$&$-1$ \\
 $\sigma_3$& $p-1$&$ p-1$&$-1$&$-1$ \\
 $\sigma_4$& $(p-1)^2$&$-( p-1)$&$-( p-1)$&$-( p-1)$ \\ \hline
\end{tabular}
\end{center}

To find the unitary matrix $U$ we use the formula written down in section 2 to obtain the $4\times4$ matrix $U$ as follows: 
$$ U=\dfrac{1}{p} \begin{bmatrix}
1&\sqrt{p-1}&\sqrt{p-1}&p-1\\\sqrt{p-1}&-1&p-1&-\sqrt{p-1}\\\sqrt{p-1}&p-1&-1&-\sqrt{p-1}\\p-1&-\sqrt{p-1}&-\sqrt{p-1}&1
\end{bmatrix}$$

At this point it is convenient to consider the general case of $ G=\mathbb{Z}_p^d $, 
\begin{center}
$\Gamma=\lbrace\begin{bmatrix}
a_1&&&&\\
&a_2&&\bf{0}&\\
&&\ddots&&\\
&\bf{0}&&\ddots&\\
&&&&a_d
\end{bmatrix}\mid a_i\in\mathbb{Z}_p^\times \rbrace$
\end{center}
 the diagonal subgroup of order $(p-1)^d$ of $GL_d(\mathbb{Z}_p)$. 

Orbits of $\Gamma$ on $G$ are as follows: 
\\
$Y_1=\lbrace(0,0,\ldots,0)\rbrace$ is one orbit. Let $y^{(k)}=(1^k,0^{d-k})$ be a vector of $G$ with $k$ one's in different positions. Then $y^{(k)}\Gamma$ consists of vectors with non-zero entries in exactly $k$ different positions. Therefor the orbit $y^{(k)}$ has size $(p-1)^k$. Since this $k$ positions is taken out of $d$ positions, therefore we have $\begin{pmatrix}
d\\k
\end{pmatrix}$ orbits of this shapes each of size $(p-1)^k$. Hence we have $\underset{k=0}{\overset{d}{\sum}}\begin{pmatrix}
d\\k
\end{pmatrix}=2^d$  orbits of $\Gamma$ on $G$. Each orbit has size $(p-1)^k$. Since $\underset{k=0}{\overset{d}{\sum}}\begin{pmatrix}
d\\k
\end{pmatrix}(p-1)^k=p^d=\vert G\vert$ , all the orbits are counted. 
  
Orbits of $\Gamma$ on $Irr(G)$ have the same setting as above. In this case if $\psi_x$ is a representative of the orbit $X$ of $\Gamma$ on $Irr(G)$, then we may assume $x=x^{(l)}=(1^l,0^{d-l})$ is a vector with $l$ ones in different positions, hence. 
$$ \sigma_X(y)=\underset{x\in X}{\sum}\psi_{x}(y)=\underset{x\in X}{\sum}e(\dfrac{x\cdot y}{p}) $$
And it is computable if the inner product $x\cdot y$ is known.

\section{J-Symmetric groups}
Let $G=\mathbb{Z}_n^d$ and $\Gamma$ be a subgroup of $GL_d(\mathbb{Z}_n)$. By \cite{BR} we have to assume that $\Gamma$ is symmetric, i.e. $\Gamma=\Gamma^t$, in order to conclude that the action of $\Gamma$ on $G$ and on $Irr(G)$ generate the same orbits. Most of the results on supercharacter theory of $G$ holds if we assume $\Gamma$ is J-Symmetric. Suppose there is a fixed symmetric invertible matrix $J\in GL_d(\mathbb{Z}_n)$ such that $J\Gamma=\Gamma^t J$. As before the action of $\Gamma$ on $G$ is by $y\longmapsto yA$ and by identifying $\psi_x\in Irr(G)$ with $x$, the action of $\Gamma$ on $Irr(G)$ is by $x\longmapsto xA^{-t}$ for $A\in \Gamma$. 
 
If $(\mathcal{X},\mathcal{Y})$ is the supercharacter theory obtained in this way, then we set $\mathcal{X}=\lbrace X_1,X_2,\ldots,X_h\rbrace$ and $\mathcal{Y}=\lbrace Y_1,Y_2,\ldots,Y_h\rbrace$ and $\sigma_i=\sigma_{X_i}$, $1\leq i\leq r$, then the unitary matrix $U$ is replaced by $U=\dfrac{1}{\sqrt{n^d}}\begin{bmatrix}\dfrac{\sigma_i(Y_j)\sqrt{\vert Y_j\vert}}{\sqrt{\vert X_i\vert}}\end{bmatrix}^h_{i,j=1}$ .

In this section we consider $G=\mathbb{Z}_p^3$, $p$ a prime, and $J=\begin{bmatrix}
0&0&1\\0&1&0\\1&0&0\
\end{bmatrix}$, 
\begin{center}
$\Gamma=\lbrace\begin{bmatrix}
a&b&c\\0&d&b\\0&0&a\
\end{bmatrix}\mid a,b\in \mathbb{Z}_p^\times, b,c\in\mathbb{Z}_p\rbrace$
\end{center}
 is a subgroup of $GL_3(\mathbb{Z}_p)$ of order $p^2(p-1)^2$. It is obvious that $\Gamma$ is a J-Symmetric group. 

Orbits of $\Gamma$ on $G$ are as follows: 
\begin{align*}
& Y_1=\lbrace(0,0,0)\rbrace\\
& Y_2=(0,0,1)\Gamma=\lbrace(0,0,a)\mid a\in \mathbb{Z}_p^\times\rbrace\\
& Y_3=(0,1,0)\Gamma=\lbrace(0,d,b)\mid d\in \mathbb{Z}_p^\times,b\in\mathbb{Z}_p\rbrace\\
& Y_4=(1,0,0)\Gamma=\lbrace(a,b,c)\mid a\in \mathbb{Z}_p^\times,b,c\in\mathbb{Z}_p\rbrace.
\end{align*}

We have 
\begin{align*}
&\vert Y_1\vert=1\\&\vert Y_2\vert=p-1\\&\vert Y_3\vert=p(p-1)\\&\vert Y_4\vert=p^2(p-1).
\end{align*} 
Since $\vert Y_1\vert+\vert Y_2\vert+\vert Y_3\vert+\vert Y_4\vert=p^3$, we deduce that $Y_1,Y_2,Y_3$ and $Y_4$ are orbits of $\Gamma$  on $G$. It is easy to see that the orbits of $\Gamma$  on $Irr(G)$ are as follows: 
\begin{align*}
& X_1=\lbrace(0,0,0)\rbrace\\
& X_2=(1,0,0)\Gamma=\lbrace(a,0,0)\mid a\in \mathbb{Z}_p^\times\rbrace\\
& X_3=(0,1,0)\Gamma=\lbrace(a,b,0)\mid a\in \mathbb{Z}_p,b\in\mathbb{Z}_p^\times\rbrace\\
& X_4=(1,0,0)\Gamma=\lbrace(a,b,c)\mid a,b\in \mathbb{Z}_p,c\in\mathbb{Z}_p^\times\rbrace.
\end{align*}
We have 
\begin{align*}
&\vert X_1\vert=1\\&\vert X_2\vert=p-1\\&\vert X_3\vert=p(p-1)\\&\vert X_4\vert=p^2(p-1).
\end{align*} 

Let the supercharacter associated to $X_i$ be $\sigma_i$. The following supercharacter table for the group $G$ is constructed: 
\begin{center}
\begin{tabular}{|c|c|c|c|c|}\hline
 $\mathbb{Z}_p^3$& $Y_1$&$ Y_2$&$Y_3$&$Y_4$ \\ 
  superclass size & $1$&$p-1$&$p(p-1)$&$p^2(p-1)$ \\ \hline
 $\sigma_1$& $1$&$ 1$&$1$&$1$ \\
 $\sigma_2$& $p-1$&$ p-1$&$p-1$&$-1$ \\
 $\sigma_3$& $p(p-1)$&$ p(p-1)$&$-p$&$0$ \\
 $\sigma_4$& $p^2(p-1)$&$-p^2$&$0$&$0$ \\ \hline
\end{tabular}
\end{center}
The unitary table associated with the above table is: 
$$ U=\dfrac{1}{p\sqrt{p}}\begin{bmatrix}
1&\sqrt{p-1}&\sqrt{p(p-1)}&p\sqrt{p-1}\\\sqrt{p-1}&p-1&(p-1)\sqrt{p}&-p\\\sqrt{p(p-1}&(p-1)\sqrt{p}&-p&0\\p\sqrt{p-1}&-p&0&0
\end{bmatrix} $$

As general case let us consider $G=\mathbb{Z}_p^d$, 
\begin{center}
$\Gamma=\lbrace\begin{bmatrix}
1&a_2&a_3&\cdots&a_d\\
0&1&a_2&\cdots&a_{d-1}\\
\vdots&&&&\vdots\\
0&0&\cdots&\cdots&a_2\\
0&0&\cdots&\cdots&1\\
\end{bmatrix}\mid a_i\in \mathbb{Z}_p, 2\leq i\leq d\rbrace$
\end{center}
which is J-symmetric with respect to the $d\times d$ matrix $J=\begin{bmatrix}
&&&&&1\\
&&\bf{0}&&1&\\
&&&.&&\\
&&.&&&\\
&.&&\bf{0}&&\\
1&&&&&
\end{bmatrix}$. We have $\vert\Gamma\vert=p^{p-1}$ and it is a $p$-group. 

The orbits of $\Gamma$ on $G=\mathbb{Z}_p^d$ are grouped as follows: 
\begin{align*}
&Y_1=\lbrace(0,0,\ldots,0)\rbrace\\
&Y_2=(\alpha,0,\ldots,0)\Gamma=\lbrace(\alpha,a_2,a_3,\ldots,a_d)\mid a_i\in \mathbb{Z}_p\rbrace,\alpha\in\mathbb{Z}_p^\times\\
\end{align*}
hence $Y_2$ is the union of $p-1$ orbits each of size $p^{d-1}$. 
\begin{align*}
Y_3=(0,\alpha,0,\ldots,0)\Gamma=\lbrace(0,\alpha,\alpha_2,\ldots,\alpha_{d-1})\mid\alpha_i\in \mathbb{Z}_p\rbrace,\alpha\in\mathbb{Z}_p^\times
\end{align*}
hence $Y_3$ is the union of $p-1$ orbits each of size $p^{d-2}$. If we continue in this way we obtain
\begin{align*}
Y_d=(0,0,\ldots,\alpha,0)\Gamma=\lbrace(0,0,\ldots,\alpha,a_{d})\mid a_2\in \mathbb{Z}_p\rbrace
\end{align*}
has size $p$ and is the union of $p-1$ orbits, 
\begin{align*}
Y_{d+1}=(0,0,\ldots,\alpha)\Gamma=\lbrace(0,0,\ldots,\alpha)\rbrace
\end{align*}
is the union of $p-1$ orbits each of size $1$. 

Since $1+(p-1)(p^{d-1}+p^{d-2}+\ldots+1)=p^d=\vert G\vert$, all the orbits of $\Gamma$ on $G$ are counted. Therefore there are $1+(p-1)d$ orbits. 

To find the shapes of the orbits of $\Gamma$ on $Irr(G)$, we mention that each irreducible character of $G$ has degree $1$. We have $Irr(G)=\lbrace \psi_x\mid x\in G\rbrace$ which may be represented by elements $x$ of $G$ under the action $x\longmapsto xA^{-t}$ where $A\in \Gamma$. Therefore we obtain the following orbits: 
\begin{align*}
&X_1=\lbrace(0,0,\ldots,0)\rbrace\\
&X_2=(\alpha,0,\ldots,0)\Gamma=\lbrace(\alpha,0,\ldots,0)\rbrace\\
&X_3=(0,\alpha,0,\ldots,0)\Gamma=\lbrace(a_2,\alpha,0,\ldots,0)\mid a_2\in\mathbb{Z}_p\rbrace\\
&\vdots\\
&X_{d+1}=(0,0,\ldots,\alpha)\Gamma=\lbrace(a_2,a_3,\ldots,a_d,\alpha)\mid \alpha_i\in\mathbb{Z}_p\rbrace.
\end{align*}
Each set $X_i$, $2\leq i\leq d+1$ is the union of $p-1$ orbit each of size $p^{i-2}$. 

Now if $\psi_x$ is a representative of the orbit $X$ of $\Gamma$ on $Irr(G)$, we may choose 
\begin{center}
$x=(0,0,\ldots,0,\ldots,0)\in X_i$
\end{center}
 hence if $\sigma_X$ is the supercharacter associated to $X$, then for $y\in Y_j$ we have 
\begin{center}
$\sigma_X(y)=\underset{x\in X}{\sum}\psi_x(y)=\underset{x\in X}{\sum}e(\dfrac{x\cdot y}{p})$
\end{center}
if $x$ and $y$ are taken from orbits such that $x\cdot y=0$, then $\underset{x\in X}{\sum}e(\dfrac{x\cdot y}{p})=\vert X\vert=p^{i-2}$, provided $X=X_i$. Otherwise if  $x\cdot y\neq0$ we obtain $\underset{x\in X}{\sum}e(\dfrac{x\cdot y}{p})=0$. In this way the supercharacter table of $G$ is computed.

\end{document}